# The General Form Of Cyclic Orthonormal Generators In $\mathbb{R}^N$

Kerry M. Soileau

June 6, 2005

ABSTRACT

In this paper we give a definition of cyclic orthonormal generators (cogs) in $\mathbb{R}^N$. We give a general canonical form for their expression. Further, we give an explicit formula for computing the canonical form of any given cog.

Each vector in $\mathbb{R}^N$ can be expressed in the form $(a_1, a_2, \cdots, a_N)$ with each $a_i$ real. (Here and in the following we assume $N$ is a whole number not less than 2). We define the cycle of a vector $(a_1, a_2, \cdots, a_N)$ to be the vector $(a_2, a_3, \cdots, a_N, a_1)$. If we begin with a vector, cycle it, cycle the result, and so on, we will eventually return to the original vector, having "generated" at most $N-1$ cyclic companions. For example the vector $(0,0,1,0,\cdots,0) \in \mathbb{R}^N$ has $N-1$ cyclic companions, while $(1,1,\cdots,1) \in \mathbb{R}^N$ has none.

The purpose of this paper is to answer the question, under what conditions does a vector in $\mathbb{R}^N$ together with its cyclic companions form an orthonormal basis for $\mathbb{R}^N$? We call such a vector a cyclic orthonormal generator in $\mathbb{R}^N$. Henceforth in this paper we will use the abbreviation cog. Note that any cyclic companion to a cog is itself a cog.

Are there any nontrivial examples? Yes; with some early trial and error the author found that $\left(\frac{2}{3}, \frac{2}{3}, -\frac{1}{3}\right)$ and its cyclic companions $\left(\frac{2}{3}, -\frac{1}{3}, \frac{2}{3}\right)$ and $\left(-\frac{1}{3}, \frac{2}{3}, \frac{2}{3}\right)$ form an orthonormal basis for $\mathbb{R}^3$. This led to a more systematic search with a computer, and features of many cogs found by computer inspired this paper.

We begin by introducing the following

Notation: For any integer n, let $(n)^*$ be the unique integer satisfying

$(n)^* \in \{1, 2, \cdots, N\}$

and

$(n)^* \equiv n \bmod N$ .

This notation permits us compactly to define the characteristic property of a cog in $\mathbb{R}^N$. The cyclic companions of such a vector $(a_1, a_2, \cdots, a_N)$ are just the vectors $\left(a_{(1+s)^*}, a_{(2+s)^*}, \cdots, a_{(N+s)^*}\right)$ where $s \in \{1, 2, \ldots, N-1\}$. By means of the summation formula for dot products, it is immediately seen that a vector $(a_1, a_2, \cdots, a_N)$ is a cog in $\mathbb{R}^N$ if and only if $\sum_{j=1}^{N} a_j a_{(j+s)^*} = \delta_0^s$ for every $s \in \{0, 1, 2, \ldots, N-1\}$.

Next, suppose $(a_1, a_2, \cdots, a_N)$ is a cog in $\mathbb{R}^N$. Then define for $n \in \{0, 1, 2, \ldots, N-1\}$,

$$C_n = \sum_{j=1}^{N} a_j \cos\left(\frac{\pi}{4} - \frac{2\pi n}{N}(j-1)\right)$$

and

$$S_n = \sum_{j=1}^{N} a_j \sin\left(\frac{\pi}{4} - \frac{2\pi n}{N}(j-1)\right)$$

Lemma 1. $(C_n, S_n)$ lies on the unit circle for each $n \in \{0, 1, \cdots, N-1\}$.

Proof: Fix $n \in \{0, 1, \cdots, N-1\}$. Then

$$C_n^2 + S_n^2 = \overline{(C_n + iS_n)}(C_n + iS_n)$$

$$=\overline{\left(\sum_{j=1}^{N} a_j \left(\cos\left(\frac{\pi}{4}-\frac{2\pi n}{N}(j-1)\right)+i\sin\left(\frac{\pi}{4}-\frac{2\pi n}{N}(j-1)\right)\right)\right)}$$

$$\cdot\left(\sum_{k=1}^{N} a_k \left(\cos\left(\frac{\pi}{4}-\frac{2\pi n}{N}(k-1)\right)+i\sin\left(\frac{\pi}{4}-\frac{2\pi n}{N}(k-1)\right)\right)\right)$$

$$=\overline{\left(\sum_{j=1}^{N} a_j \exp\left(\left(\frac{\pi}{4}-\frac{2\pi n}{N}(j-1)\right)i\right)\right)}\cdot\left(\sum_{k=1}^{N} a_k \exp\left(\left(\frac{\pi}{4}-\frac{2\pi n}{N}(k-1)\right)i\right)\right)$$

$$=\left(\sum_{j=1}^{N} a_j \exp\left(\left(-\frac{\pi}{4}+\frac{2\pi n}{N}(j-1)\right)i\right)\right)\cdot\left(\sum_{k=1}^{N} a_k \exp\left(\left(\frac{\pi}{4}-\frac{2\pi n}{N}(k-1)\right)i\right)\right)$$

$$\sum_{j=1}^{N}\sum_{k=1}^{N} a_j a_k \exp\left(\frac{2\pi n}{N}(k-j)i\right)$$

$$=\sum_{s=0}^{N-1}\sum_{j=1}^{N}\sum_{\substack{k=1\\ k-j\equiv s \bmod N}}^{N} a_j a_k \exp\left(\frac{2\pi n}{N}(k-j)i\right)$$

$$=\sum_{s=0}^{N-1}\sum_{j=1}^{N}\sum_{\substack{k=1\\ k\equiv j+s \bmod N}}^{N} a_j a_k \exp\left(\frac{2\pi n}{N}si\right)=\sum_{s=0}^{N-1}\sum_{j=1}^{N} a_j a_{(j+s)^*} \exp\left(\frac{2\pi n}{N}si\right)$$

$$=\sum_{s=0}^{N-1}\exp\left(\frac{2\pi ns}{N}i\right)\sum_{j=1}^{N} a_j a_{(j+s)^*}$$

which by our assumption that $(a_1, a_2, \cdots, a_N)$ is a cog in $\mathbb{R}^N$,

$$=\sum_{s=0}^{N-1}\exp\left(\frac{2\pi ns}{N}i\right)\delta_0^s = \exp(0) = 1.$$

Thus $C_n^2 + S_n^2 = 1$, and thus $(C_n, S_n)$ lies on the unit circle.

<u>Note</u>: The following proofs will employ the formulas

$$\sum_{n=0}^{N-1}\cos\left(\frac{2\pi na}{N}+\theta\right)=\begin{cases} N\cos\theta & (N|a) \\ 0 & \neg(N|a)\end{cases}$$

and

$$\sum_{n=0}^{N-1}\sin\left(\frac{2\pi na}{N}+\theta\right)=\begin{cases}N\sin\theta & (N|a)\\ 0 & \neg(N|a)\end{cases}$$

where $a$ is a nonzero integer. These are obtained as special cases of equations 361.9 and 361.8, respectively, in $[1]$.

<u>Theorem 1</u>: If $(a_1,a_2,\cdots,a_N)$ is a cog in $\mathbb{R}^N$, then there exist $\{\theta_n\}_{n=0}^{N-1}$ such that for $k\in\{1,2,\cdots,N\}$,

$$a_k=\frac{\sqrt{2}}{N}\sum_{n=0}^{N-1}\cos\left(\frac{2\pi n}{N}(k-1)+\theta_n\right).$$

<u>Proof</u>: Suppose $(a_1,a_2,\cdots,a_N)$ is a cog in $\mathbb{R}^N$, and choose $\{\theta_n\}_{n=0}^{N-1}$ such that for $n\in\{0,1,\cdots,N-1\}$

$$\cos\theta_n=\sum_{j=1}^{N}a_j\cos\left(\frac{\pi}{4}-\frac{2\pi n}{N}(j-1)\right)$$

$$\sin\theta_n=\sum_{j=1}^{N}a_j\sin\left(\frac{\pi}{4}-\frac{2\pi n}{N}(j-1)\right).$$

Lemma 1 guarantees the possibility of such choices. Fix $k\in\{1,2,\cdots,N\}$. Then

$$\frac{\sqrt{2}}{N}\sum_{n=0}^{N-1}\cos\left(\frac{2\pi n}{N}(k-1)+\theta_n\right)$$
$$=\frac{\sqrt{2}}{N}\sum_{n=0}^{N-1}\cos\left(\frac{2\pi n}{N}(k-1)\right)\cos\theta_n-\frac{\sqrt{2}}{N}\sum_{n=0}^{N-1}\sin\left(\frac{2\pi n}{N}(k-1)\right)\sin\theta_n$$

$$=\frac{\sqrt{2}}{N}\sum_{n=0}^{N-1}\cos\left(\frac{2\pi n}{N}(k-1)\right)\sum_{j=1}^{N}a_j\cos\left(\frac{\pi}{4}-\frac{2\pi n}{N}(j-1)\right)$$
$$-\frac{\sqrt{2}}{N}\sum_{n=0}^{N-1}\sin\left(\frac{2\pi n}{N}(k-1)\right)\sum_{j=1}^{N}a_j\sin\left(\frac{\pi}{4}-\frac{2\pi n}{N}(j-1)\right)$$

$$=\frac{\sqrt{2}}{N}\sum_{j=1}^{N}a_j\sum_{n=0}^{N-1}\cos\left(\frac{2\pi n}{N}(k-1)+\frac{\pi}{4}-\frac{2\pi n}{N}(j-1)\right)$$

$$=\frac{\sqrt{2}}{N}\sum_{j=1}^{N}a_j\sum_{n=0}^{N-1}\cos\left(\frac{2\pi n}{N}(k-j)+\frac{\pi}{4}\right)$$

which by formula and the observation that $N\,|(k-j)$ only if $k=j$,

$$=\frac{\sqrt{2}}{N}\sum_{j=1}^{N}a_j\left(N\cos\frac{\pi}{4}\right)\delta_j^k=a_k$$, as desired.

Thus, we have shown that every cog $(a_1,a_2,\cdots,a_N)$ satisfies

$$a_k=\frac{\sqrt{2}}{N}\sum_{n=0}^{N-1}\cos\left(\frac{2\pi n}{N}(k-1)+\theta_n\right)$$

for some $\{\theta_n\}_{n=0}^{N-1}$.

At this point a natural question is, what conditions must the $\{\theta_n\}_{n=0}^{N-1}$ satisfy? An interesting property is demonstrated by the following

Lemma 2: If $(a_1,a_2,\cdots,a_N)$ is a cog in $\mathbb{R}^N$, and $\{\theta_n\}_{n=0}^{N-1}$ are such that for $n\in\{0,1,\cdots,N-1\}$,

$$\cos\theta_n=\sum_{j=1}^{N}a_j\cos\left(\frac{\pi}{4}-\frac{2\pi n}{N}(j-1)\right)$$

and

$$\sin\theta_n=\sum_{j=1}^{N}a_j\sin\left(\frac{\pi}{4}-\frac{2\pi n}{N}(j-1)\right),$$

then

$$\theta_0\equiv\begin{cases}\frac{\pi}{4}\bmod 2\pi & \left(\sum_{j=1}^{N}a_j=1\right)\\ \frac{5\pi}{4}\bmod 2\pi & \left(\sum_{j=1}^{N}a_j=-1\right)\end{cases}$$

and

$$\theta_n + \theta_{N-n} \equiv \frac{\pi}{2} \bmod 2\pi \text{ for } n \in \{1, 2, \cdots, N-1\}.$$

<u>Proof</u>: Fix $n \in \{1, 2, \cdots, N-1\}$. Then note that

$$\cos(\theta_n + \theta_{N-n}) = \cos\theta_n \cos\theta_{N-n} - \sin\theta_n \sin\theta_{N-n}$$
$$= \sum_{j=1}^{N} a_j \cos\left(\frac{\pi}{4} - \frac{2\pi n}{N}(j-1)\right) \sum_{k=1}^{N} a_k \cos\left(\frac{\pi}{4} - \frac{2\pi(N-n)}{N}(k-1)\right)$$
$$- \sum_{j=1}^{N} a_j \sin\left(\frac{\pi}{4} - \frac{2\pi n}{N}(j-1)\right) \sum_{k=1}^{N} a_k \sin\left(\frac{\pi}{4} - \frac{2\pi(N-n)}{N}(k-1)\right)$$

$$= \sum_{j=1}^{N} \sum_{k=1}^{N} a_j a_k \cos\left(\frac{\pi}{2} - \frac{2\pi n}{N}(j-k) - 2\pi(k-1)\right)$$

which since $k-1$ ranges over whole numbers,

$$= \sum_{j=1}^{N} \sum_{k=1}^{N} a_j a_k \cos\left(\frac{\pi}{2} - \frac{2\pi n}{N}(j-k)\right) = \sum_{j=1}^{N} \sum_{k=1}^{N} a_j a_k \sin\left(\frac{2\pi n}{N}(j-k)\right)$$
$$= \sum_{s=0}^{N-1} \sum_{k=1}^{N} a_{(k+s)^*} a_k \sin\frac{2\pi n s}{N}$$

$$= \sum_{s=0}^{N-1} \sin\frac{2\pi n s}{N} \sum_{k=1}^{N} a_{(k+s)^*} a_k = \sum_{s=0}^{N-1} \sin\frac{2\pi n s}{N} \delta_0^s = 0.$$

Thus

$$\cos(\theta_n + \theta_{N-n}) = 0.$$

Next, note that

$$\sin(\theta_n + \theta_{N-n}) = \sin\theta_n \cos\theta_{N-n} + \cos\theta_n \sin\theta_{N-n}$$

$$= \sum_{j=1}^{N} a_j \sin\left(\frac{\pi}{4} - \frac{2\pi n}{N}(j-1)\right) \sum_{k=1}^{N} a_k \cos\left(\frac{\pi}{4} - \frac{2\pi(N-n)}{N}(k-1)\right)$$
$$+ \sum_{j=1}^{N} a_j \cos\left(\frac{\pi}{4} - \frac{2\pi n}{N}(j-1)\right) \sum_{k=1}^{N} a_k \sin\left(\frac{\pi}{4} - \frac{2\pi(N-n)}{N}(k-1)\right)$$

$$= \sum_{j=1}^{N} \sum_{k=1}^{N} a_j a_k \sin\left( \frac{\pi}{2} - \frac{2\pi n}{N} (j-k) - 2\pi (k-1) \right)$$

which since $k-1$ ranges over whole numbers,

$$= \sum_{j=1}^{N} \sum_{k=1}^{N} a_j a_k \sin\left( \frac{\pi}{2} - \frac{2\pi n}{N} (j-k) \right) = \sum_{j=1}^{N} \sum_{k=1}^{N} a_j a_k \cos\left( \frac{2\pi n}{N} (j-k) \right)$$

$$= \sum_{s=0}^{N-1} \sum_{k=1}^{N} a_{(k+s)^*} a_k \cos \frac{2\pi ns}{N}$$

$$= \sum_{s=0}^{N-1} \cos \frac{2\pi ns}{N} \sum_{k=1}^{N} a_{(k+s)^*} a_k = \sum_{s=0}^{N-1} \cos \frac{2\pi ns}{N} \delta_0^s = \cos 0 = 1 .$$

Thus

$$\sin(\theta_n + \theta_{N-n}) = 1 .$$

Since $\cos(\theta_n + \theta_{N-n}) = 0$, we infer immediately that $\theta_n + \theta_{N-n} \equiv \frac{\pi}{2} \bmod 2\pi$, as desired.

To complete the proof, we now demonstrate the possible values of $\theta_0$.

First note that

$$\cos\theta_0 = \sum_{j=1}^{N} a_j \cos \frac{\pi}{4} = \frac{1}{\sqrt{2}} \sum_{j=1}^{N} a_j$$

and

$$\sin\theta_0 = \sum_{j=1}^{N} a_j \sin \frac{\pi}{4} = \frac{1}{\sqrt{2}} \sum_{j=1}^{N} a_j .$$

We see immediately that $\cos\theta_0 = \sin\theta_0$. This implies

$$\theta_0 \equiv \begin{cases} \dfrac{\pi}{4} \bmod 2\pi & \left(\displaystyle\sum_{j=1}^{N} a_j = 1\right) \\ \dfrac{5\pi}{4} \bmod 2\pi & \left(\displaystyle\sum_{j=1}^{N} a_j = -1\right) \end{cases}$$

and the proof is complete.

Corollary: If $(a_1, a_2, \cdots, a_N)$ is a cog in $\mathbb{R}^N$, then $\sum_{j=1}^{N} a_j = \pm 1$.

We have so far shown that any cog in $\mathbb{R}^N$ can be expressed in a canonical form (as in the statement of Theorem 1) using $\{\theta_n\}_{n=0}^{N-1}$ such that $\theta_0 \equiv \frac{\pi}{4}$ or $\frac{5\pi}{4} \bmod 2\pi$ and $\theta_n + \theta_{N-n} \equiv \frac{\pi}{2} \bmod 2\pi$ for $n \in \{1, 2, \cdots, N-1\}$.

Example: Recall the cog $\left(\frac{2}{3}, \frac{2}{3}, -\frac{1}{3}\right)$. We now use Theorem 1 to put this cog into canonical form. To do this we must find appropriate values for $\theta_0$, $\theta_1$, and $\theta_2$. In this example we have $a_1 = \frac{2}{3}$, $a_2 = \frac{2}{3}$, and $a_3 = -\frac{1}{3}$. This implies $\sum_{j=1}^{N} a_j = 1$, hence $\theta_0 = \frac{\pi}{4}$. Next,

$$\cos\theta_1 = \frac{2}{3}\cos\left(\frac{\pi}{4}\right) + \frac{2}{3}\cos\left(-\frac{5\pi}{12}\right) - \frac{1}{3}\cos\left(-\frac{13\pi}{12}\right) \approx 0.965926$$

$\sin\theta_1 = \frac{2}{3}\sin\left(\frac{\pi}{4}\right) + \frac{2}{3}\sin\left(-\frac{5\pi}{12}\right) - \frac{1}{3}\sin\left(-\frac{13\pi}{12}\right) \approx \text{-}0.258819$ hence $\theta_1 \approx 6.021386$. Next,

$$\cos\theta_2 = \frac{2}{3}\cos\left(\frac{\pi}{4}\right) + \frac{2}{3}\cos\left(-\frac{13\pi}{12}\right) - \frac{1}{3}\cos\left(-\frac{29\pi}{12}\right) \approx -.258819$$

$\sin\theta_2 = \frac{2}{3}\sin\left(\frac{\pi}{4}\right) + \frac{2}{3}\sin\left(-\frac{13\pi}{12}\right) - \frac{1}{3}\sin\left(-\frac{29\pi}{12}\right) \approx .965926$ hence $\theta_2 \approx 1.832596$.

Theorem 2: If we define

$$a_k = \frac{\sqrt{2}}{N}\sum_{n=0}^{N-1}\cos\left(\frac{2\pi n}{N}(k-1)+\theta_n\right) \text{ for } k\in\{1,2,\cdots,N\},$$

where $\theta_0 \equiv \frac{\pi}{4}$ or $\frac{5\pi}{4}\bmod 2\pi$ and $\theta_n+\theta_{N-n}\equiv\frac{\pi}{2}\bmod 2\pi$ for $n\in\{1,2,\cdots,N-1\}$,

then $(a_1,a_2,\cdots,a_N)$ is a cog in $\mathbb{R}^N$.

<u>Proof</u>: Define $a_k = \frac{\sqrt{2}}{N}\sum_{n=0}^{N-1}\cos\left(\frac{2\pi n}{N}(k-1)+\theta_n\right)$ for $k\in\{1,2,\cdots,N\}$, with $\{\theta_n\}_{n=0}^{N-1}$ chosen so as to satisfy the conditions $\theta_0 \equiv \frac{\pi}{4}$ or $\frac{5\pi}{4}\bmod 2\pi$ and $\theta_n+\theta_{N-n}\equiv\frac{\pi}{2}\bmod 2\pi$ for $n\in\{1,2,\cdots,N-1\}$.

Now fix $s\in\{0,1,\cdots,N-1\}$. Observe that

$$\sum_{j=1}^{N} a_j a_{(j+s)^*}$$

$$=\sum_{j=1}^{N}\frac{\sqrt{2}}{N}\sum_{n=0}^{N-1}\cos\left(\frac{2\pi n}{N}(j-1)+\theta_n\right)\frac{\sqrt{2}}{N}\sum_{m=0}^{N-1}\cos\left(\frac{2\pi m}{N}\left((j+s)^*-1\right)+\theta_m\right)$$

$$=\frac{2}{N^2}\sum_{n=0}^{N-1}\sum_{m=0}^{N-1}\sum_{j=1}^{N}\frac{1}{2}\left\{\begin{array}{l}\cos\left(\frac{2\pi n}{N}(j-1)+\theta_n+\frac{2\pi m}{N}\left((j+s)^*-1\right)+\theta_m\right)\\+\cos\left(\frac{2\pi n}{N}(j-1)+\theta_n-\frac{2\pi m}{N}\left((j+s)^*-1\right)-\theta_m\right)\end{array}\right\}$$

which since $(j+s)^*\equiv(j+s)\bmod N$ and $m$ is an integer,

$$=\frac{2}{N^2}\sum_{n=0}^{N-1}\sum_{m=0}^{N-1}\sum_{j=1}^{N}\frac{1}{2}\left\{\begin{array}{l}\cos\left(\frac{2\pi n}{N}(j-1)+\theta_n+\frac{2\pi m}{N}(j-1+s)+\theta_m\right)\\+\cos\left(\frac{2\pi n}{N}(j-1)+\theta_n-\frac{2\pi m}{N}(j-1+s)-\theta_m\right)\end{array}\right\}$$

$$=\frac{1}{N^2}\sum_{n=0}^{N-1}\sum_{m=0}^{N-1}\sum_{j=1}^{N}\left\{\begin{array}{l}\cos\left(\frac{2\pi}{N}(n+m)(j-1)+\frac{2\pi m}{N}s+\theta_n+\theta_m\right)\\+\cos\left(\frac{2\pi}{N}(n-m)(j-1)-\frac{2\pi m}{N}s+\theta_n-\theta_m\right)\end{array}\right\}$$

$$=\frac{1}{N^2}\sum_{n=0}^{N-1}\sum_{m=0}^{N-1}\sum_{j=0}^{N-1}\left\{\begin{array}{l}\cos\left(\frac{2\pi}{N}(n+m)j+\frac{2\pi m}{N}s+\theta_n+\theta_m\right)\\+\cos\left(\frac{2\pi}{N}(n-m)j-\frac{2\pi m}{N}s+\theta_n-\theta_m\right)\end{array}\right\}.$$

Before continuing with the computation we make several observations based upon the summation formulae quoted earlier from [1].

The internal sum

$$\sum_{j=0}^{N-1}\cos\left(\frac{2\pi}{N}(n+m)j+\frac{2\pi m}{N}s+\theta_n+\theta_m\right) \text{ vanishes unless } N\,\big|\,(n+m),$$

i.e., $n=m=0$ or $n+m=N$. In case $n=m=0$, the sum is easily seen to be $N\cos 2\theta_0=0$ by our assumption on $\theta_0$. In case $n+m=N$, the sum is

$$\sum_{j=0}^{N-1}\cos\left(2\pi j+\frac{2\pi ms}{N}+\theta_{N-m}+\theta_m\right)=\sum_{j=0}^{N-1}\cos\left(2\pi j+\frac{2\pi ms}{N}+\frac{\pi}{2}\right)$$

since $j$ takes on integer values and $\theta_n+\theta_m\equiv\frac{\pi}{2}$ by assumption.

This sum reduces further to

$$\sum_{j=0}^{N-1}-\sin\left(2\pi j+\frac{2\pi ms}{N}\right)=-\sum_{j=0}^{N-1}\sin\frac{2\pi ms}{N}=-N\sin\frac{2\pi ms}{N}.$$

Thus

$$\sum_{j=0}^{N-1}\cos\left(\frac{2\pi}{N}(n+m)j+\frac{2\pi m}{N}s+\theta_n+\theta_m\right)=\begin{cases}-N\sin\frac{2\pi ms}{N} & (n+m=N)\\ 0 & (n+m\neq N)\end{cases}$$

Next, the internal sum

$$\sum_{j=0}^{N-1}\cos\left(\frac{2\pi}{N}(n-m)j-\frac{2\pi ms}{N}+\theta_n-\theta_m\right) \text{ vanishes unless } N\,\big|\,(n-m),$$

i.e., $n=m$. In this case the sum reduces to

$N\cos\frac{2\pi ms}{N}$.

Thus

$$\sum_{j=0}^{N-1}\cos\left(\frac{2\pi}{N}(n-m)j-\frac{2\pi ms}{N}+\theta_n-\theta_m\right)=\begin{cases}N\cos\frac{2\pi ms}{N} & (n=m)\\ 0 & (n\neq m)\end{cases}$$

Recall that

$$\sum_{j=1}^{N}a_j a_{(j+s)*}=\frac{1}{N^2}\sum_{n=0}^{N-1}\sum_{m=0}^{N-1}\sum_{j=0}^{N-1}\left\{\begin{array}{l}\cos\left(\frac{2\pi}{N}(n+m)j+\frac{2\pi m}{N}s+\theta_n+\theta_m\right)\\ +\cos\left(\frac{2\pi}{N}(n-m)j-\frac{2\pi m}{N}s+\theta_n-\theta_m\right)\end{array}\right\}$$

$$=\frac{1}{N^2}\sum_{n=0}^{N-1}\sum_{\substack{m=0\\ n+m=N}}^{N-1}\left(-N\sin\frac{2\pi ms}{N}\right)+\frac{1}{N^2}\sum_{n=0}^{N-1}\sum_{\substack{m=0\\ n=m}}^{N-1}\left(N\cos\frac{2\pi ms}{N}\right)$$

$$=\frac{1}{N^2}\sum_{n=0}^{N-1}-N\sin\frac{2\pi(N-n)s}{N}+\frac{1}{N^2}\sum_{n=0}^{N-1}N\cos\frac{2\pi ns}{N}$$

$$=-\frac{1}{N}\sum_{n=0}^{N-1}\sin\frac{2\pi(N-n)s}{N}+\frac{1}{N}\sum_{n=0}^{N-1}\cos\frac{2\pi ns}{N}$$

$$=-\frac{1}{N}\sum_{n=0}^{N-1}-\sin\frac{2\pi ns}{N}+\frac{1}{N}\sum_{n=0}^{N-1}\cos\frac{2\pi ns}{N}$$

$$=\frac{1}{N}\sum_{n=0}^{N-1}\sin\frac{2\pi ns}{N}+\frac{1}{N}\sum_{n=0}^{N-1}\cos\frac{2\pi ns}{N}=\frac{1}{N}0+\frac{1}{N}N\delta_0^S=\delta_0^S$$

Thus $\sum_{j=1}^{N}a_j a_{(j+s)*}=\delta_0^s$. Since we fixed an arbitrary $s\in\{0,1,\cdots,N-1\}$, by definition $(a_1,a_2,\cdots,a_N)$ is a cog in $\mathbb{R}^N$.

<u>Definition</u>: If $(a_1,a_2,\cdots,a_N)$ is a cog in $\mathbb{R}^N$, and $\{\theta_n\}_{n=0}^{N-1}$ is such that

$$a_k = \frac{\sqrt{2}}{N} \sum_{n=0}^{N-1} \cos\left( \frac{2\pi n}{N}(k-1) + \theta_n \right) \text{ for } k \in \{1, 2, \cdots, N\},$$

then we call $\{\theta_n\}_{n=0}^{N-1}$ a <u>representation</u> of $(a_1, a_2, \cdots, a_N)$.

Naturally we ask: does a given cog in $\mathbb{R}^N$ have more than one representation? The answer is found in

<u>Theorem 3</u>: If $\{\theta_n\}_{n=0}^{N-1}$ and $\{\tau_n\}_{n=0}^{N-1}$ are representations of a cog in $\mathbb{R}^N$, then $\theta_n \equiv \tau_n \bmod 2\pi$ for $n \in \{0, 1, \cdots, N-1\}$. Conversely, if $\{\theta_n\}_{n=0}^{N-1}$ is a representation of a cog in $\mathbb{R}^N$, and $\theta_n \equiv \tau_n \bmod 2\pi$ for $n \in \{0, 1, \cdots, N-1\}$, then $\{\tau_n\}_{n=0}^{N-1}$ is also a representation of the same cog.

<u>Proof</u>: Suppose $(a_1, a_2, \cdots, a_N)$ is a cog in $\mathbb{R}^N$, and suppose that $\{\theta_n\}_{n=0}^{N-1}$ and $\{\tau_n\}_{n=0}^{N-1}$ are representations of $(a_1, a_2, \cdots, a_N)$ in $\mathbb{R}^N$. Fix $k \in \{1, 2, \cdots, N\}$. Since $\{\theta_n\}_{n=0}^{N-1}$ is a representation of $(a_1, a_2, \cdots, a_N)$ in $\mathbb{R}^N$, we have by definition

$$a_k = \frac{\sqrt{2}}{N} \sum_{n=0}^{N-1} \cos\left( \frac{2\pi n}{N}(k-1) + \theta_n \right).$$

Fix $m \in \{0, 1, \cdots, N-1\}$. Then

$$\cos\theta_m = \sum_{j=1}^{N} a_j \cos\left( \frac{\pi}{4} - \frac{2\pi m}{N}(j-1) \right)$$

and

$$\sin\theta_m = \sum_{j=1}^{N} a_j \sin\left( \frac{\pi}{4} - \frac{2\pi m}{N}(j-1) \right).$$

Since $\{\tau_n\}_{n=0}^{N-1}$ is a representation of $(a_1, a_2, \cdots, a_N)$ in $\mathbb{R}^N$, we have by definition

$$a_k = \frac{\sqrt{2}}{N}\sum_{n=0}^{N-1}\cos\left(\frac{2\pi n}{N}(k-1)+\tau_n\right).$$

Then

$$\cos\tau_m = \sum_{j=1}^{N} a_j \cos\left(\frac{\pi}{4}-\frac{2\pi m}{N}(j-1)\right)$$

and

$$\sin\tau_m = \sum_{j=1}^{N} a_j \sin\left(\frac{\pi}{4}-\frac{2\pi m}{N}(j-1)\right).$$

Hence it is immediately seen that

$$\cos\theta_m = \cos\tau_m$$

and

$$\sin\theta_m = \sin\tau_m,$$

Since we fixed an arbitrary $m \in \{0,1,\cdots,N-1\}$, it follows that

$$\theta_n \equiv \tau_n \bmod 2\pi \text{ for } n \in \{0,1,\cdots,N-1\}.$$

For the converse, suppose $\{\theta_n\}_{n=0}^{N-1}$ is a representation of $(a_1,a_2,\cdots,a_N)$ in $\mathbb{R}^N$, and $\theta_n \equiv \tau_n \bmod 2\pi$ for $n \in \{1,2,\cdots,N-1\}$. Fix $k \in \{1,2,\cdots,N\}$. Then

$$a_k = \frac{\sqrt{2}}{N}\sum_{n=0}^{N-1}\cos\left(\frac{2\pi n}{N}(k-1)+\theta_n\right).$$

Fix $m \in \{0,1,\cdots,N-1\}$. Since by hypothesis $\theta_m \equiv \tau_m \bmod 2\pi$, we have immediately that

$$\cos\left(\frac{2\pi m}{N}(k-1)+\theta_m\right) = \cos\left(\frac{2\pi m}{N}(k-1)+\tau_m\right).$$

Thus

$$a_k = \frac{\sqrt{2}}{N}\sum_{n=0}^{N-1}\cos\left(\frac{2\pi n}{N}(k-1)+\theta_n\right) = \frac{\sqrt{2}}{N}\sum_{n=0}^{N-1}\cos\left(\frac{2\pi n}{N}(k-1)+\tau_n\right)$$

and the proof is complete.

CONCLUSIONS

In this paper we give a definition of cyclic orthonormal generators (cogs) in $\mathbb{R}^N$. We give a general canonical form for their expression. Further, we give an explicit formula for computing the canonical form of any given cog. We remark in closing that for a given choice of $N \geq 3$, the set of cogs in $\mathbb{R}^N$ is homeomorphic with the $\lfloor \frac{N-1}{2} \rfloor$-manifold $\left\{\frac{\pi}{4}, \frac{5\pi}{4}\right\} \times [0, 2\pi)^{\lfloor \frac{N-1}{2} \rfloor}$. This manifold can be regarded as the union of two non-intersecting $\lfloor \frac{N-1}{2} \rfloor$-torii.

REFERENCES

1. Adams, E. P., and Hippisley, R. L., Smithsonian Mathematical Formulae and Tables of Elliptic Functions, Smithsonian Institution, Washington, D. C., 1922.